\documentclass[letterpaper, 10 pt, conference]{ieeeconf}  
\IEEEoverridecommandlockouts                              
\usepackage{mathptmx} 
\usepackage{times} 
\usepackage{amsmath} 
\usepackage{amssymb}  
\usepackage{cite}
\title{\LARGE \bf
Large Gain Stability of Extremum Seeking\\
and Higher-Order Averaging Theory 
}

\usepackage{stmaryrd}
\usepackage{amsopn}
\usepackage{dsfont}
\usepackage{graphicx}
\usepackage{color}
\usepackage{array}
\usepackage[T1]{fontenc}
\usepackage[utf8]{inputenc}
\author{Gabriel \textsc{Bousquet} and  Jean-Jacques \textsc{Slotine}%
\thanks{Nonlinear Systems Laboratory, MIT, Cambridge, USA. This work was supported  in part by a Fulbright Science and Technology fellowship.}%
}

\usepackage{lipsum}
\usepackage{ntheorem}

\newtheorem{lemma}{Lemma}

\newtheorem*{fundamentalDecomposition}{Fundamental Decomposition}

\begin{document}
\maketitle

{\abstract
Convergence of Extremum Seeking (ES) algorithms has been established in the limit of small gains.
Using averaging theory and contraction analysis, we  propose a framework for computing explicit bounds on the departure of the ES scheme from its ideal dominant-order average dynamics. The bounds remain valid for possibly large gains. They allow us to establish stability and estimate convergence rates, and they open the way to selecting "optimal" finite  gains for the ES scheme.
}

\section{Introduction}

Extremum Seeking (ES) is a special class of algorithms designed to optimize dynamic systems\cite{Tan2010}. Typically, traditional unconstrained optimization algorithms are designed to find the maximum of a map $h:\mathbf{R}^n\to\mathbf{R}$ through a sequence of evaluations $h(x_n)$.
In the kind of problem addressed by ES, the map parameters are continuous functions of time $x(t)\in\mathbf{R}$ and the map output is continuously measured, maybe through a dynamic system $h([z], x)$ where $[z]$ represents the possible internal dynamics\cite{ariyur2003}.

Extremum Seeking can be traced back to\cite{Leblanc1922} and was an active field until the third quarter of the 20th century\cite{Sternby1980}. While dormant during the following decades, it gained much attention since the early 2000's, in part after the theoretical progresses in\cite{Krstic2000,Tan2006}, where  the first proof of  stability for the ES scheme is provided. ES has known a rapid development in the last decade and its range of applications is expanding rapidly with the generalization of large scale and low cost autonomous dynamic systems and robots\cite{Leyva2006,zhang2012}.

Design and stability analyses of ES schemes have been presented for several classes of systems that follow a similar structure\cite{Nesic2010}:
\begin{itemize}
\item a particular gradient related method is to be mimicked. The gradient and perhaps the Hessian\cite{Ghaffari2012} of the system at its operating point are estimated  by introducing oscillatory perturbations to the input and measuring the correlation with the evolution of the output and are used to define a direction of search, 
\item The loop is closed by setting operating point $x$ to drift on a slow timescale along the search direction.
\item By design, the slow time dynamics of $x$ is a minimization direction so that the slow time equation typically reduces to $d_t\overline{x} = k\nabla h(\overline{x})$ or $kH^{-1}\nabla h$ (the over-line represents the time average over a dither period).
\item From there, averaging theory and singular perturbation are invoked and it is shown that if the gains used in the design of the ES scheme are small enough, $x$ converges to a limit cycle in neighborhood of the optimal point $x^*$\cite{Tan2006}.
\end{itemize}
By doing so, only the stability in the limit of infinitely small gains is proven. The speed of convergence however is driven by the amplitude of those small gains. For practical applications it is therefore important to know how large they  can be while convergence is maintained, search precision is quantified and a satisfactory search speed is provided. Although is it important to mention that some basic work and scaling remarks have been done in this direction \cite{Rotea2000,Tan2006,Nesic2010,Teel2001,Krstic2000a}, no explicit solution has been provided yet, let alone for a nonlinear analysis.

The objective of this paper is to show that those limitations can be overcome if two quantitative tools --averaging theory up to higher orders and singular serturbation theory revisited by contraction theory\cite{Lohmiller1999,Chen2012,DelVecchio2013}-- are used. The basic idea is to bound or estimate the departure of the real system from the ideal, designed system, and to build an auxiliary optimization problem called meta-optimization that will select the small gains so that the errors remain below a fixed error. Just as the small gain stability has been proven for specific algorithms (as opposed to the whole class of ES algorithms), the same goes for the finite gain theory. In part 2 we summarize the averaging theory as presented in Chapter 3 of \cite{Sanders} and discuss its application to finite order averaging before stating the results from the contraction analysis of singularly perturbed systems\cite{DelVecchio2013} and showing how they can be combined for the analysis of ES schemes. In part 3 we apply those quantitative tools to the one dimension map maximization, first in its most simple form as in section 4.1 in \cite{Tan2006}, and then in its high and low pass filtered version \cite{Krstic2000}. We then extend the first case to the optimization of a black box dynamic system. We finish part 3 by illustrating on the simple 2-dimensional case how higher order averaging can be used as a qualitative design tool to avoid high order, undesired terms.

We conclude this introduction by mentioning that the analysis presented here for the sake of simplicity on 1-D objective maps extends straightforwardly to multidimensional cases, should apply seamlessly to applications such as the non-holonomic search\cite{Cochran2009} and constitute an interesting direction of research in the more general setup of stochastic ES\cite{Liu2010}.
\section{Theory}

\subsection{Averaging}

\subsubsection{Fundamental relation of averaging}

Let $f$, $w$ be smooth vector fields in $\mathbf{R}^n$. We note
$\mathcal{D}_{w} h = \nabla h\cdot w$ and
  $\mathcal{L}_wh = \mathcal{D}_wh
  - \mathcal{D}_hw$. The results from averaging theory (\cite{Sanders} Chapter 3 and  \cite{Murdock}, Annex C) that are used in the present paper rely on the following relations: consider the dynamic system
\begin{equation}
\dot{x} = f(x)\label{eq:nonlineq}
\end{equation}
and the new variable $y$ defined by the implicit transformation
\begin{equation}
x = e^{\mathcal{D}_w}y\label{eq:averaging0}
\end{equation}
The dynamics of $y$ given by equation $\dot{y}=g_\infty(y)$ with
\begin{equation}
 g_\infty(y) = e^{\mathcal{L}_w}{f(y)}\label{eq:averaging1}
\end{equation}
where the exponential of operator has its usual definition:
\[x = y + \mathcal{D}_{w} y 
  + \frac{(\mathcal{D}_{w})^2y}{2!}
  + \frac{(\mathcal{D}_{w})^3y}{3!}
+ \dots
\]
\[
g_\infty(y) = f(y) + \mathcal{L}_wf(y) 
+ \frac{(\mathcal{L}_{w})^2f(y)}{2!}
+ \frac{(\mathcal{L}_{w})^3f(y)}{3!}
+ \dots
\]
Simply put, this theorem means that given a variable change $x = U_\infty(y) = e^{\mathcal{D}_w}y$ (the introduction of the generator function $w$ can be thought of as an indirect way to define $U_\infty$), the dynamics for the new variable takes the simple form given by equation~(\ref{eq:averaging1}).

\subsubsection{Interpretation when $f$ and $w$ are written as series}

\label{sec:series_exp1}

If ${f}({x})= \sum_{i=1}^\infty
\epsilon^i{f}_i({x})$ and ${w}({y})=
\sum_{i=1}^\infty \epsilon^i{w}_i({y})$, where $\epsilon$ is a given
scalar, the previous relations must hold for all orders in $\epsilon$.
Collecting the terms of same order together, the previous expressions become
\[
x = U_\infty(y) = y + \sum_{i=1}^\infty \epsilon^i u_i(y)
\]
and
\[
g_\infty(y) = \sum_{i=  1}^\infty\epsilon^i g_i(y)
\]
The $u_i$'s are sums and products of $w_k$'s up to order $i$ and their derivatives up to order $i-1$ while the $g_i$'s are sums and products of $f_j$'s up to order $i$, $w_j$'s up to order $i-1$ and their derivatives up to order $i-1$.
In the present paper, we don't venture deeper into the properties of the $g_i$'s and $w_i$'s. At this point it is important to  mention that, although their expressions are non trivial (especially in the case of non-autonomous systems, which are of our interest), they can be computed algorithmically at each order.
\subsubsection{Application to a non autonomous system}

Let's consider the system
\begin{equation}
  \label{eq:non_autonom_1}
\dot{{x}} = {f(x}, t)  
\end{equation}
and augment it with $x_{n+1} = t$ into its autonomous form:
\[ {X} =
\begin{bmatrix}
  {x}\\
  t
\end{bmatrix}
\qquad
\dot{{X}} = {F(X)} =
\begin{bmatrix}
  {f(x}, t)\\
  1
\end{bmatrix}
\]
The relations of the previous section still hold, and the differential operators $\mathds{L}$ and $\mathds{D}$ of the $n+1$ dimensional, spatiotemporal system can be rewritten as functions of their  $n$ dimensional, spatial counterparts.
Assuming  ${W} = [{w}^T \ 0]^T$ we get
\[
{(\mathds{D}_{W}})^p{F} =
\begin{bmatrix}
  \mathcal{D}_{w}^p{f}\\
  0
\end{bmatrix} \qquad p\geq 1
\]
\[
{(\mathds{L}_{W}}) {F}=
\begin{bmatrix}
  \mathcal{L}_{w}{f} - \partial_t {w}\\
  0
\end{bmatrix}
\]
\[
{(\mathds{L}_{W}})^p {F}=
\begin{bmatrix}
  \mathcal{L}_{w}\\
  0
\end{bmatrix}{(\mathds{L}_{W}})^{p-1} {F }\qquad n \geq 2
\]

Finally, if we note $\tilde{\mathcal{L}}$ the modified non autonomous, Lie
bracket such that $\tilde{\mathcal{L}} = \mathcal{L}_{w}(\cdot)
- \partial_t{w}$ and $\tilde{\mathcal{L}}^p =
\mathcal{L}\tilde{\mathcal{L}}^{p-1}$, $p\geq 2$, the non autonomous system
can be transformed by
\[{x} = e^{\mathcal{D}_{w}}{y}= U_\infty(y,t) = y + \delta U_\infty(y, t) \]
into
\[\dot{{y}} = {g_\infty(y}, t)\]
with
\begin{equation}
  {g_\infty(y}, t) = e^{\tilde{\mathcal{L}}_{w}}{f(y}, t)\label{eq:averaging2}
\end{equation}

\subsubsection{$f$ and $w$ as series expansions, take 2}
Let's assume that ${f(x}, t)$ and ${w(x},
t)$ can be written as series in $\epsilon$, as in
part~\ref{sec:series_exp1}. It is noticeable that, in the right hand
side of
equation~(\ref{eq:averaging2}), $\partial_t{w}_i$ appears first
at order $i$ in the Taylor expansion while ${w}_i$ doesn't appear until
order $\epsilon^{i+1}$. At each order in $\epsilon$ in
equation~(\ref{eq:averaging2}), we therefore have:
\[{g}_i({y}, t) = -\partial_t{w}_i +
{E}_i({y}, t)
\]
where ${E}_i({y}, t)$ can be expressed as sums, products and
derivatives of ${f}_j$, $j\le i$,  ${w}_j$, $j\le i-1$ and their derivatives.
In other words, by choosing
\begin{equation}
  \label{eq:algo1}
  {w}_i({y}, t) = \int_0^t {E}_i({y}, t) dt + {K}_i({y)}
\end{equation}
${g}_i$ is independent of time at each order. At this point, it
is important to note that if ${f}$ is $T-$periodic (as will be
assumed in the following),  so are the
${E}_i$'s, ${w}_i$'s and, as a consequence,
${U_\infty}$. Also, if the ${K}_i$'s are chosen to be ${0}$,
${y}$ and ${x}$ coincide at $ t = 0 [T]$. A probably more appropriate choice, used in the present article, is to set $K_i$ so that $\overline{u}_i = 0$ for $i\ge 1$. It leads to $\overline{x} = y$.

From this, we can conclude that there exists an algorithmic way to transform the non autonomous system
$\dot{{x}} = \sum_{i\geq 1} \epsilon^i {f}_i({x}, t)$
into an autonomous $\dot{{y}} = \sum_{i\ge 1} \epsilon^i {g}_i({y})$ at each order in $\epsilon$ by choosing an appropriate change of variable. That change of variables is periodic in time if $f$ is.


\subsubsection{Finite order averaging}

The previous sections provide an algorithmic way to build a change of variables that transforms the non autonomous, periodically driven, dynamic system into an autonomous one. For practical applications however, it is important to note that the variable change can only be performed up to a finite order $n$.
Consider that the previous algorithm  has been carried out up to order $n$ so that $w_i$ $i\le n$ have been defined. We note $\tilde{w}$ the truncation of $w$. Then, $e^{\mathcal{D}_{\tilde{w}}}y$ can be computed and truncated into $U_n(y, t)= U(y, t)$ that defines the relationship between $y$ and $x$. This change of variables is then used to compute the dynamics of $y$. Derivating
\[
x = U(y, t)
\]
we get
\[
f(U(y, t), t) = \partial_yU\dot{y} + \partial_tU
\]
Assuming that $U$ remains invertible at all times:
\begin{equation}
\dot{y} = [\partial_yU]^{-1}(f(U(y, t), t) - \partial_tU) \label{eq:ydotWithU}
\end{equation}
By construction, the dynamics of $y$ is autonomous up to order $n$. The above equation can therefore be rewritten as
\[
\dot{y}= g(y) + R_g(y, t)
\]
where $g$ is the truncation of $g_\infty$ and $R_g(y, t) = O(\epsilon^{n+1})$ arises from the fact that the averaging is only carried up to order $n$.

\subsection{Singular Perturbation with Contraction}

\subsubsection{Contraction theory\cite{Lohmiller1999}}
Consider the system $\dot{x} = f(x, t)$. It is said to be contracting if all trajectories converge exponentially towards each other. A necessary condition for contraction is that there exists a metric $\Theta(x, t)$ such that $\Theta^T\Theta$ is uniformly positive definite and $\beta>0$ such that
\[
F = \dot{\Theta}\Theta^{-1} + \Theta \nabla f \Theta^{-1} \preccurlyeq -\beta I
\]

$\beta$ is called the contraction rate. We also define $\chi$ as a bound on the condition number of $\Theta$. A useful lemma from contraction theory is the robustness lemma
\begin{lemma}
  if $f$ is contracting with rate $\beta$ and $R$ is a bounded perturbation such that $\dot{y} = f(y) + R(t)$, 
  $y(t)$ converges to a $|R|/\kappa$ neighborhood of $x(t)$ where $\kappa = \beta/\chi$. We say that the system is $\kappa$-robust.
\end{lemma}

\subsubsection{Singular perturbation\cite{DelVecchio2013}}

Consider the dynamic system
  \begin{align*}
    \nu \dot{z}&= g(x, z)\\
    \dot{x} &= f(x, z, t)\\
  \end{align*}
  \begin{lemma}
    Assume that the fast system $\nu \dot{z}= g(x_0(t), z)$ is
    partially $\lambda/\nu$-robust with respect to $z$ (the $x-z$
    coupling has been replaced by an external forcing). Write
    $\gamma(x)$ its equilibrium, assume that there exists $d>0$ such
    that $|\partial_x\gamma(x) f(x, z, t)| \le d$. Assume also
    that $f$ is Lipschitz in $z$ with constant $\alpha$. Assume
    eventually that $z(t = 0) = \gamma(x(t= 0))$ and let $x_\gamma$ be
    the solution of the reduced singular perturbation system
    $\dot{x}_\gamma = f(x_\gamma, \gamma(x_\gamma), t)$ ($y$ is frozen
    to its equilibrium state). Then,
\[
|\dot{x} - \dot{x}_\gamma|\le \frac{d\alpha\nu}{\lambda_z}
\]
\end{lemma}

\subsection{Utilizing higher order averaging  and singular perturbations with contraction in Extremum Seeking}

The typical ES scheme is as follows. Let's denote $x$ the parameters to be optimized, $\xi$ the internal dynamics built in ES (filters/estimation) and $z$ the internal dynamics of the unknown system. The state equation is then:
\begin{equation}
  \label{eq:fullES}
  \frac{d}{dt}
  \begin{bmatrix}
    z\\
    x\\
    \xi
  \end{bmatrix}
  =
  \begin{bmatrix}
    \frac{1}{\nu}A(z, x)\\
    \epsilon B_{\epsilon, r}\left(C(z), \xi, t)\right)\\
    \epsilon D_{\epsilon, r}\left(C(z), \xi, t)\right)
  \end{bmatrix}
\end{equation}
where $B, D$ are nonlinear functions dependent on the particular ES system. While ES usually depends on several small parameters (gains and cutoff frequencies), we assume that those parameters have been rescaled so that the system depends on a parameter of possibly small amplitude $\epsilon$ and a set of parameters $r = [r_1\dots r_p]$ of order unity.
The slow system is
\begin{equation}
\label{eq:slowES}
  \frac{d}{dt}
  \begin{bmatrix}
    x_s\\
    \xi_s
  \end{bmatrix}
  =\epsilon
  \begin{bmatrix}
    B_{\epsilon, r}\left(h(x_s), \xi_s, t)\right)\\
    D_{\epsilon, r}\left(h(x_s), \xi_s, t)\right)
  \end{bmatrix} = f(X_s)
\end{equation}
with $X = [x; \xi]$ and $X_s = [x_s; \xi_s]$. $B$ and $D$ known, as they are the design of ES and $h = C(\gamma(\cdot)))$ is the static output of the fast dynamics black box system to be optimized. Only $h$ is unknown.

If the hypothesis of the previous sections hold, applying lemma~2 gives:
\begin{equation}
  |\dot{X} - \dot{X}_s| \le \frac{\alpha d\nu}{\lambda_z}\label{eq:errDyn}
\end{equation}
This result is useful as it allows to give a bound on the difference between a system where a map if optimized (the limiting case for which ES is designed) from when a black box dynamic system is optimized, as a function of its bandwidth.

Let's now assume that an ES scheme has been designed to operate on a map. With averaging theory, the change of variable $U$ is constructed for system~(\ref{eq:slowES}) algorithmically as explained in the previous part. Applied to system~(\ref{eq:fullES}) with the use of equations (\ref{eq:ydotWithU}) and (\ref{eq:errDyn}) this change of variables gives:
\begin{fundamentalDecomposition}
  \begin{align}
    X &= Y + \delta U(Y, t)\\
    \dot{Y} &= g_0(Y) + \delta g_s(Y) + R_{g_s}(Y, t) + R_\nu(Y, z, t)
  \end{align}
\end{fundamentalDecomposition}
$Y = [x_{av}, \xi_{av}]$, 
where $\delta U = O(\epsilon)$ is the non identity part of the coordinate transformation,  $g_0$ gathers the dominant order, ideal terms (typically, the gradient or Newton descent), $\delta g_s= o_\epsilon(g_0)$ represents the computed departure from the ideal descent that is autonomous and polynomial in $\epsilon$ up to order $n$, $R_{g_s} = O(\epsilon^{n+1})$ is the higher order error of the slow system and $R_\nu$ is the error induced by the black box dynamics.
In particular, if bounds are known on $f$ and its derivatives (gathered in the notation $\llbracket f\rrbracket= (\|f\|, \|f'\|\dots)$), each of those higher terms can be bounded: $|\delta U|\le K_1(\epsilon, r, \llbracket f\rrbracket )$, $|\delta g_s|\le K_2(\epsilon, r, \llbracket f\rrbracket )$, $|R_s|\le K_3(\epsilon, r, \llbracket f\rrbracket )$, $|R_\nu| =| \partial_yU^{-1}(\dot{X} - \dot{X}_s)| \le K_4(\nu, \epsilon, r, \llbracket f\rrbracket )$.

Using this decomposition, a relatively simple optimization can be run to select $\nu$, $\epsilon$ and the order unity parameters $r_i$ to maximize the search speed while keeping the error terms below some bounds. Particular possible meta-optimizations and useful relaxations are presented in the next section.

\section{Applications}

\subsection{1 state, 1-D system}

To illustrate this technique, we apply the procedure to the most simple 1-state ES scheme from\cite{Tan2006} and work it through step by step.
For ease of computation, we assume the dither signal to be sinusoidal, although the same analysis can be carried out with any periodic signal.
The state equation for this ES scheme is:
\begin{equation}
  \label{eq:base_eq}
  \dot{x} = -\eta h(x + a\sin t)\sin t
\end{equation}
Qualitatively, if $\eta$ is small, the right hand side is small so that $x$ is varying slowly and its long term evolution is given by the short time average of the right hand side, which reduces to $\overline{x} = -a\eta h'(\overline{x})/2$ if $a$ is also small. The system therefore mimics a  gradient descent scheme of $h$ at rate $a\eta/2$.
Here, $\eta$ and $a$ are parameters set up by the practitioner. In the cited literature, they were assumed to be small parameters, enough so that the approximations that make $x$ approximately driven by $\nabla h$ hold.
It is important however to note that the descent speed it driven by $a\eta$, and that those parameters should therefore be as large as possible so that the search performs promptly. It appears that the amplitude of $a$ and $\eta$ constitute a trade-off between speed and accuracy.
Our interest in to use averaging theory to quantify the departure from the ideal gradient descent that the system mimics when the parameters are small and compute what finite gains are acceptable to make the search perform at a set precision, while maximizing the search speed.

Before moving further, we take advantage of this example to discuss the parametrization of ES. Averaging theory, as presented in the previous section and in the reference literature, is parametrized with one small parameter only. In the case of ES, there are several small parameters so that the averaging analysis cannot be performed as is.
One way to bypass this theoretical limitation is to parametrize each parameter as a function of the averaging parameter $\epsilon$: $\eta = \eta(\epsilon)$ and $a = a(\epsilon)$. We require $a$ and $\eta$ to be of class $\mathcal{K}_\infty$. The most natural choice is to use power laws, $a= \epsilon^n$ and $\eta = p\epsilon^m$, $m, n \in \mathbf{N}^*$. This way, we separate in $\eta$ the ``magnitude'' part $\epsilon^m$, and its fine tuning value $p$. Then, to the reparametrized system
\[
  \dot{x} = -\eta(\epsilon) h\left(x + a(\epsilon)\sin t\right)\sin t
\]
corresponds an averaged system in coordinate $y$ defined by $x = U(y)$  such that
\[\dot{y} = \sum_{i= 1}^n\epsilon^i g_i(y) + R_g(y, t)
\]
where the $g_i$ are sums and products of $h$, $\eta$ and $a$ up the $i^\text{th}$ derivative.

If $m=n=1$, the averaged system to the dominant order is
\begin{align*}
\dot{y} &=-\frac{a\eta}{2}h'(y) -  \frac{1}{16}\left( a\eta^3 \mathcal{L}^2_hh' + \eta a^3 h^{(3)} \right)  + R_g \\
x & = y + \eta \sin(t)h(y) + R_u
\end{align*}
where $R_g = O(\epsilon^4)$ and $R_u = O(\epsilon^2)$ are higher order terms in $a$, $\eta$ that can be computed explicitly.
The term $\frac{a\eta}{2}h'$ is the ideal gradient descent that the scheme is intended to reproduce. The middle term in the equation for $y$ was noted $\delta g_s$ in the previous section. For a given function $h$, the speed or search is proportional to $a\eta$. Optimizing for the gains implies therefore some sort of maximization of $a\eta$. 

The dynamics for $y$ is in  form of the general dynamic nonlinear system with noise $\delta g_s + R_g$ considered in the contraction theory lemma.
Assuming that $h$ is $\kappa$-robust, lemma~1 can be applied to show that $y$ follows the ideal trajectory $\dot{z} = \frac{a\eta}{2}h'(z)$ with an error that is at most $\delta_1 = 2\frac{\|R_g\| +\|\delta g_s\|}{a\eta\kappa}$. 
Similarly, $|x - y|\le \delta_2 = \eta \|h'\| + \|R_u\|$. Assuming that an estimate on $\llbracket h\rrbracket$ (and therefore on $R_i$) in known, the gains can be chosen to maximize the search speed while keeping $\delta_1$ and $\delta_2$ below some tolerance error. In practice, it is possible to relax the constraints slightly by only considering the dominant error terms. There are several possible strategies and some are listed below:
\begin{enumerate}
\item The guaranteed meta-optimization consists in bounding the distance from the real system to the ideal system. A simple triangular inequality shows that $|x - z|\le \delta_1 + \delta_2$ so that the problem reduces to:
  \begin{align*}
    \max_{\eta, a\ge 0} &\ \eta a\\
    \text{s.t.}&\  \sum\delta_i \le \Delta
  \end{align*}
This is the safest optimization, as it guarantees convergence and the errors bounds ($x$ will converge to a $\Delta$-neighborhood of $x^*$). However, the higher order remainders are often complex so that, even if they can be bounded provided that bounds $\llbracket h\rrbracket $ are known, the bound is unlikely to be tight, leading to overly conservative regimes. Also, the error due to the average system and the error due to the oscillation of the real system may not be of the same importance for the user, so that it might be fruitful to separate them in the constraints. Those two drawbacks motivate the next two strategies
\item Splitting up the error from the average system (the DC error) $\delta_1$ and the oscillatory error $\delta_2$, the meta-optimization becomes:
  \begin{align*}
    \max_{\eta, a\ge 0} &\ \eta a\\
    \text{s.t.}&\  \delta_i \le \Delta_i
  \end{align*}
\item Neglecting the highest order of the remainder,greatly simplifies the expressions for the bounds. For instance, in the case where $\eta = p\epsilon$ and $a = \epsilon$, the 1-D problem reduces to:
  \begin{align*}
    \max_{\eta, a\ge 0} &\ \eta a\\
    \text{s.t.}&\  \frac{2}{a\eta\kappa}\left\{\frac{1}{16}(\eta a^3 \|h^{(3)}\| + a\eta^3 \|\mathcal{L}^2_hh'\|)\right\}\le \Delta_1 \\
    \text{s.t.}&\ \eta \|h\| \le \Delta_2
  \end{align*}
  In that case, since a truncation is made, a more aggressive search is recommended by the meta-optimization. Formal guarantee of convergence is lost though so that it is wise to check, after such an meta-optimization, that the neglected terms are indeed negligible.
\item The previous strategies are well suited when the objective is to locate the optimal point (such as when performing source/pollutant tracking for instance). In other situations it might be important to keep the real system  $x+a\sin t$ close to $x^*$ at all times. The extension of strategy 1) is
  \begin{align*}
    \max_{\eta, a\ge 0} &\ \eta a\\
    \text{s.t.}&\ a + \frac{2\|R_g\|}{a\eta\kappa} + \|R_u\| \le \Delta
  \end{align*}
\end{enumerate}
\emph{Remarks:}
\begin{itemize}
\item In the case $a = \epsilon$ and $\eta = p \epsilon$, the terms in $\eta a^3$ and $a \eta^3$ in the constraint for $y$ have the same order. However, if $\eta = p\epsilon^m$, $a = \epsilon^n$ and $m\neq n$, the two terms won't necessarily have the same order. Therefore, for more general $m,n$, the symmetry may be broken leading to a problem with monomial constraints
  \begin{align*}
    \max_{\eta, a\ge 0} &\ \eta a\\
    \text{s.t.}&\  \eta^{p_1}a^{q_1} \le K_1\\
    \text{s.t.}&\ \eta^{p_2} a^{q_2} \le K_2
  \end{align*}
  It has a finite solution if $p_1/q_1 < 1 < p_2/q_2$ which is $a = K_2^\frac{p_1}{q_2p_1 - q_1p_2}K_1^\frac{-p_2}{q_2p_1 - q_1p_2}$ and $\eta = K_2^\frac{-q_1}{q_2p_1 - q_1p_2}K_1^\frac{q_2}{q_2p_1 - q_1p_2}$.
\item Assume this monomial form for $a$ and $\eta$ and  $m>n$. then, in the first constraint from case 3), the first term of the first constraint is dominant. Assuming that the expansion remains the same as in the $m=n=1$ and that both constraints are active 
 brings $\eta = \Delta_2/\|h\|$ and 
 $a =  \sqrt{\frac{8\Delta_1\kappa}{\|h^{(3)}\|}}$
\item Two consistency checks should be performed to ensure meaningful results:
  \begin{itemize}
  \item are the higher order terms indeed small?
  \item is $p$ close to unity? Otherwise, $m, n$ chosen to perform the expansion may not be the appropriate ones.
  \end{itemize}
  \item Obviously, if some rough estimates on $\llbracket h \rrbracket$ are known, they can be used to compute tighter bounds on $\delta g$ and $R_i$.

\end{itemize}

\subsection{One state case, fully worked example}

Let's consider the optimization of $h(x) = - \cos(x) + x^3/6$, starting from $x = 1$. It is chosen as a toy example for all of its derivatives are bounded by 1. The cubic term is added to break the third order symmetry at the minimum $x^* = 0$. Of course, $0$ is not a global minimum,  but for $a$, $\eta$ small enough, it is still the attractor.

We solve here  the meta-optimization problem 3). The constraints require us to have estimates for the  bounds on $h$, $h'$, $h''$, $h^{(3)}$ and $\kappa$. We take $\|h^{(i)}\| \hat{=} \max_{x\in[-1, 1]} |h^{(i)}(x)|$ and $\kappa = h''(0)$. We also set $\Delta_i= 0.01$ and start with $\eta = p\epsilon$, $a = \epsilon$. It gives the numerical solution $\eta = 0.01$, $a = 0.207$ Therefore $p \approx 0.05\ll 1$. To bring $p$ closer to unity, we repeat the averaging (it has to be done up to order 7) with $\eta = p\epsilon^3$, $a = \epsilon$. The constraint for $y$ simplifies to $\frac{a^2\|h^{(3)}\|}{8\kappa}\le \Delta_1$. This time, as stated in the remarks, the problem can be solved analytically:
  \begin{align*}
    a&=  \sqrt{\frac{8\Delta_1\kappa}{\|h^{(3)}\|}}\\
    \eta&=\frac{\Delta_2}{\|h\|}
  \end{align*}
which brings with the numeric values discussed before $\eta = 0.01$ and $a = 0.209$ (and $ p = 1.09$). The reason why the result hasn't changed much is because the term in $a\eta^3$ found when $m=n=1$ in the first constraint is small already.
\begin{figure}
  \centering
  \includegraphics[width=.5\textwidth]{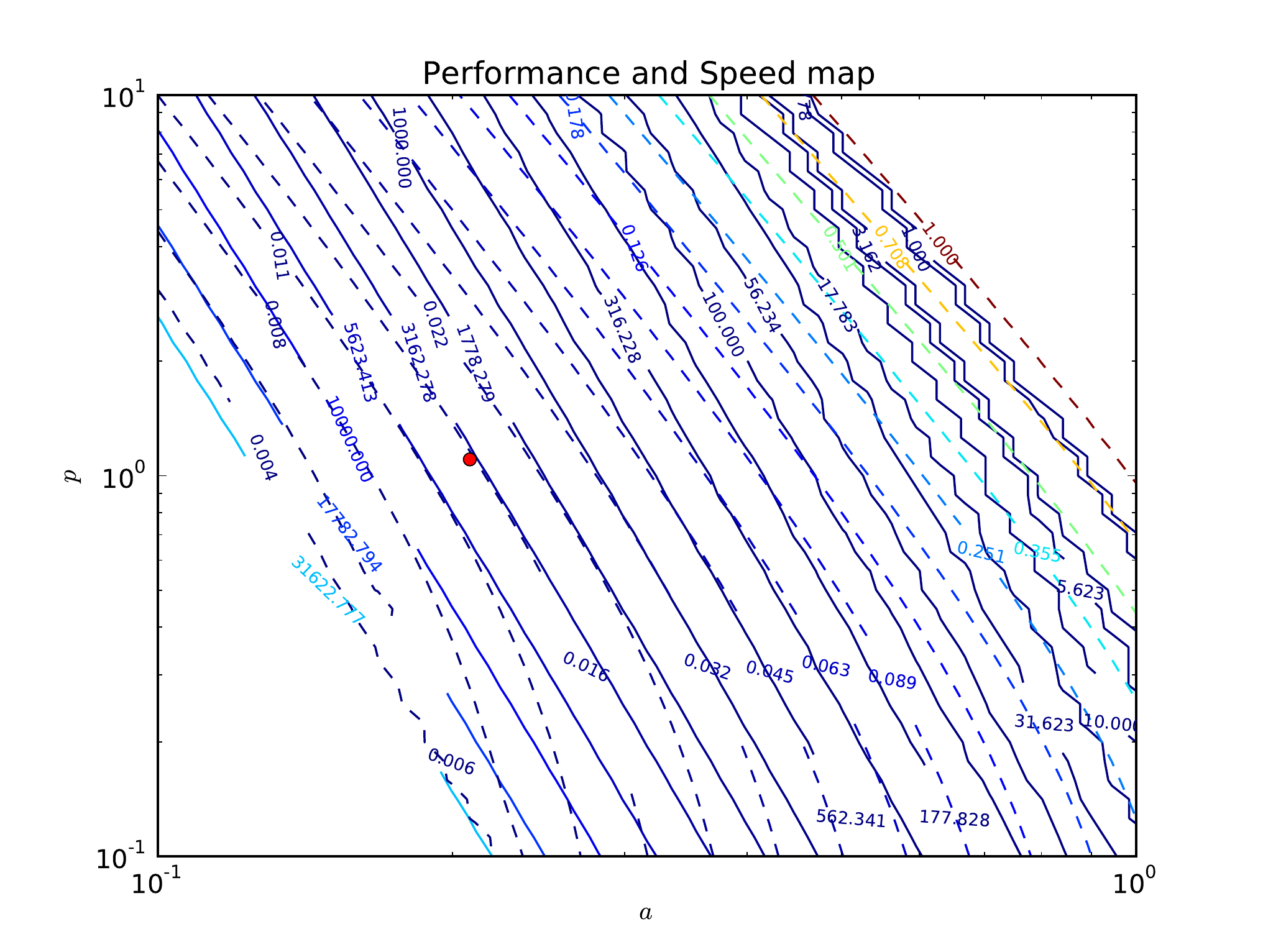}
  \caption{Performance  map of ES as a function of $a$ and $p$ (log-log scale). The full lines indicate the speed of search, the dashed lines represent the the error $|x - x^*|$. The red dot is the result from the meta optimization.}
  \label{fig:1dopt}
\end{figure}
The result is illustrated in figure~\ref{fig:1dopt}. The operation point suggested by the procedure is close to but different from the optimal point for this very map, which is expected.

\subsection{Optimization of a dynamic system}

In this part we show how the singular perturbation theory allows to set the frequency of the dither. For simplicity, we will consider again the 1-state ES, extended with a dynamic system:
\begin{equation}
  \label{eq:blackBox}
  \frac{d}{dt}
  \begin{bmatrix}
    z\\
    x
  \end{bmatrix}
  =
  \begin{bmatrix}
    -z + x \\
    k h(z)\sin\omega t
  \end{bmatrix}
\end{equation}
In the present article we illustrate the technique with a first order filter for ease of computation  but insist that the technique applies to any nonlinear, contracting system. To put system (\ref{eq:blackBox}) in form~(\ref{eq:fullES}), we change the timescale $\tau = \omega t$ and introduce $y = x + a\sin\tau$ and the reduced parameter $\eta = k/\omega$:
\begin{equation}
  \frac{d}{d\tau}
  \begin{bmatrix}
    z\\
    y
  \end{bmatrix}
  =
  \begin{bmatrix}
    \frac{1}{\omega}(-z + y) \\
    \eta h(z)\sin\tau + a\cos \tau
  \end{bmatrix}
\end{equation}
Lemma 2 can be applied with $d = \|h\|\eta + a$, $\lambda_z = 1$, $\nu = \omega$, $\alpha = \|h'\|$, to get
\[|d_\tau y - d_\tau y_s|\le \eta\omega\|h'\|(\eta \|h\| + a)\]
which can be transformed back into the $x, t$ variables:
\[|\dot{x} - \dot{x_s}|\le \eta\omega^2\|h'\|(\eta \|h\| + a)\]
We form the meta-optimization problem by relaxing $\|h'\|$ into $|h'(x)|$. This gives an approximation for the worst case speed at the dominant order
\[
|g_0| - |\dot{x} - \dot{x_s}| \approx \left(\frac{a\eta\omega}{2} - \eta\omega^2(\eta \|h\| + a)\right)|h'|
\]
which defines the objective. Keeping the constraints on $\delta_i$ from the previous part unchanged, the meta-optimization problem can be written as:
  \begin{align*}
    \max_{\eta, a, \omega\ge 0} &\ \left(\frac{a\eta\omega}{2} - \eta\omega^2(\eta |h| + a)\right)\\
    \text{s.t.}&\  \frac{1}{8\kappa}a^2 \|h^{(3)}\|\le \Delta_1 \\
    \text{s.t.}&\ \eta \|h\| \le \Delta_2
  \end{align*}
which gives $\omega = 0.48$ and keeps $a$ and $\eta$ unchanged.

\subsection{1-D map optimization with filtering}

Consider the map optimization with first order low and high pass filters adapted from\cite{Krstic2000}:
\begin{equation}\label{lphp}
  \begin{split}
    \dot{x} &= -\eta j\\
    \dot{\hat{h}} &= \mu\left(h(x+au) - \hat{h}\right)\\
    \dot{j} &= \gamma\left\{- \frac{a}{2}j +\left(h(x+au) - \hat{h}\right)u \right\}\\
  \end{split}
\end{equation}
with the dither $u = \sin t$ where $x$ is the parameter to be optimized, $\hat{h}$ is the estimate of $h(x)$ and $j$ is the estimate of $h'(x)$ in the ES algorithm. Assuming that the parameters $a$, $\eta$, $\mu$, $\gamma$ are all of order $\epsilon$, the average system becomes
\begin{align}
  \label{eq:avglphp}
  \begin{split}
    \dot{\tilde{h}}_{av} &= \dot{h}_{av}-\mu\left\{\tilde{h}_{av} + \frac{a^2}{4} h''\right\} + O(\epsilon^3)\\
    \dot{\tilde{j}}_{av} &= -\frac{a\gamma}{2}\tilde{j}_{av} +\eta h''j_{av} - \frac{\eta\gamma^2}{2}\tilde{h}h' \\
    &\quad +\frac{a\gamma}{2}\left\{-\mu^2h' + \eta\mu jh'' + \frac{a^2}{8}h^{(3)}\right\} 
      + O(\epsilon^5)\\
    \dot{x}_{av} &= -\eta j_{av} +O(\epsilon^4)\\
  \end{split}
\end{align}
\[x = x_{av} + \eta\gamma\tilde{h}\sin t + O(\epsilon^4)\]
where we have introduced $\tilde{h} = \hat{h} - h(x_{av})$ and $\tilde{j} = j - h'(x_{av})$.

Even more so that in the previous examples, there are several ways to optimize for the ES scheme. One is to maximize $\eta$ while keeping the dominant orders of $x - x_{av}$  and $\tilde{j}/\kappa$ below some fixed value $\Delta_1$ and $\Delta_2$. This is done  by applying Lemma 1 do the dynamics of $\tilde{h}$ and $\tilde{j}$. Taking the same $h$ as in the previous part and estimating the bounds on $|h|$\dots$|h^{(3)}|$ with their true maximum on $[-1, 1]$, we get $(a, \eta, \mu, \gamma) = (0.33, 8.8\cdot 10^{-3}, 0.093, 3.8)$. Numerical experiments show that the system converges into the set error bound in about 100 oscillations.

\subsection{ES in two dimensions as an illustration of channels interaction}

Practical cases of problems of optimization are multi-dimensional. The added complexity in multi-dimensional problems arises from two phenomena:
\begin{itemize}
\item the complexity related to ``higher dimension optimization'', independently from the dynamic system/extremum seeking nature of the problem
\item the coupling between directions in the estimation of the derivatives of $h$
\end{itemize}
In this section, we show that higher order averaging brings an elegant criterion for choosing the dithers that limit the coupling between channels. The result provided here is a 2-particle coupling similar to the 3-particle coupling from\cite{Rotea2000,Ghaffari2012}. For simplicity, we limit ourselves to the bidimensional case where $h(x_1, x_2)$ is to be minimized:
\begin{align*}
  \dot{x}_1 &= \eta d_1h(x_1 + ad_1, x_2 + ad_2)\\
  \dot{x}_2 &= \eta d_2h(x_1 + ad_1, x_2 + ad_2)\\
\end{align*}
We start with a qualitative presentation. In order to estimate the derivative of $h$ with respect to each of the variables, the correlation between the dither signal $d_i$ and the output $h$ is measured. 
For instance he average equation for $x_1$ is:
\[
\dot{\overline{x}}_1 \approx a\eta \left\{  \overline{d_1d_1}\partial_{x_1}h +  \overline{d_1  d_2}\partial_{x_2}h\right\}
\]
This qualitative analysis suggests that as long as $\overline{d_1 d_2} = 0$, the dynamics should reduce to the desired one:
\[ \langle \dot{x}_i\rangle \approx a\eta\langle d_i^2\rangle\partial_{x_i}h\]
This advocates for the possibility to use both sines and cosines for the dither, which is \emph{a priori} good as it would allow a bandwidth twice smaller for a given number of channels. However, the averaged equations for $d_1 = \cos t$ and $d_2 = \sin t$ are:
\begin{align*}
  \dot{y}_1 &= \frac{a\eta}{2}\partial_{x_1}h - \frac{\eta^2}{2}\partial_{x_2} h + O(\epsilon^3)\\
  \dot{y}_2 &= \frac{a\eta}{2}\partial_{x_2}h +  \frac{\eta^2}{2}\partial_{x_1} h + O(\epsilon^3)\\
\end{align*}
Although the term in $a\eta$ is the expected gradient, the term in $\eta^2$ is a precessing term, which, if dominant, keeps the system moving on level sets instead of following the gradient. The existence of this term advocates against the use of both sines and cosines in the ES (even if the system to be optimized has no dynamics), as it produces an undesired coupling between the channels at a dominant order.

\section{Conclusion}

In this paper, we have shown how contraction theory applied to singular perturbation and modern averaging theory can help bring qualitative and quantitative insights into Extremum Seeking. In particular we have shown methods for selecting the finite gains of ES schemes optimally. Although most techniques were presented for the one dimensional objective functions, they extend to $n$-dimensional problems. Further work includes extending the study to formal $n$ and to stochastic ES. The authors also believe that it constitutes a fruitful framework for designing new efficient ES schemes, particularly with adaptive gains. Lastly, it constitutes a basis for comparing optimally tuned ES schemes with other optimization techniques.

\bibliography{large_gain_article.bib}
\bibliographystyle{IEEEtran}

\end{document}